\documentclass[12pt]{article}
\usepackage{amsthm,amssymb}

\newtheorem{theorem}{Theorem}

\def\card#1{\vert #1 \vert}
\def\gpindex#1#2{\card {#1\colon #2}}
\def\irr#1{{\rm  Irr}(#1)}

\def\cent#1#2{{\bf C}_{#1}(#2)}

\def\ker#1{{\rm ker} (#1)}

\def\phi{{\varphi}}
\def\epsilon{{\varepsilon}}

\begin{document}



\title{A group with three real irreducible characters:
       answering a question of Moret\'o and Navarro}

\author {
       Mark L.\ Lewis
    \\ {\it Department of Mathematical Sciences, Kent State University}
    \\ {\it Kent, Ohio 44242}
    \\ E-mail: lewis@math.kent.edu
       }
\date{October 2, 2008}

\maketitle

\begin{abstract}
In this paper, we construct a group with three real irreducible
characters whose Sylow $2$-subgroup is an iterated central extension
of a Suzuki $2$-group.  This answers a question raised by Moret\'o
and Navarro who asked whether such a group exists.

MSC primary: 20C15
\end{abstract}


\section{Introduction}

In their note \cite{threereal}, Moret\'o and Navarro study groups
$G$ that have $3$ real-valued irreducible characters.  They prove
that $G$ is solvable, and it either has a cyclic Sylow $2$-subgroup
or it has a normal Sylow $2$-subgroup that is either homocyclic,
quaternion of order $8$, or an iterated central extension of a
Suzuki $2$-group whose center is an elementary abelian $2$-group.

Recall that $P$ is a Suzuki $2$-group if $P$ is a nonabelian
$2$-group, $P$ has more that one involution, and $P$ has a solvable
group of automorphisms that transitively permutes the involutions in
$P$.  Suzuki $2$-groups are studied extensively in Section VIII.7 of
\cite{HBII}.  The group $P$ is an iterated central extension of a
Suzuki $2$-group, if $P/ Z_n (P)$ is a Suzuki $2$-group for some $n
\ge 2$ (where $Z_n (P)$ is the $n$th term in the upper central
series for $P$).

Now, they present the following examples in that paper.  First,
$S_3$ is an example of such a group with a nonnormal, cyclic Sylow
$2$-subgroup. Also, ${\rm SL} (2,3)$ is such a group with the
quaternions as a normal Sylow $2$-subgroup.  They also produce an
example with a normal Sylow $2$-subgroup that is a Suzuki $2$-groups
(not of type A).  Finally, they include an example with a normal
Sylow $2$-subgroup that is a central extension of a Suzuki
$2$-group.

Moret\'o and Navarro did not have an example where the Sylow
$2$-subgroup $P$ had $P/Z_n (P)$ is a Suzuki $2$-group for $n \ge 2$
where $Z_n (P)$ is the $n$th term in the upper central series, and
they ask whether the word ``iterated'' can be removed from the
statement of their theorem. In this note, we will show that this
cannot be done by presenting an example of a group $G$ with three
real valued characters whose Sylow $2$-subgroup $P$ satisfies $P/Z_2
(P)$ is a Suzuki $2$-group.  At this time, we have not determined
whether or not there can be any examples with $n \ge 3$.

We will first discuss the construction used to obtain our group. We
then prove that this group has exactly three real irreducible
characters.  In an appendix, we include the MAGMA code for PC-group
arising from our construction, and we discuss the character table of
this group.

\section{The construction}

We now work to construct a group $G$ with $3$ real valued
irreducible characters whose Sylow $2$-subgroup is an iterated
Suzuki $2$-group.  We begin with a group first introduced by Isaacs
in \cite{coprime}.  Let $F$ be the field of order $8$.

We first establish some facts regarding $F$.  We define the {\it
trace map} $Tr : F \rightarrow Z_2$ by $Tr (\alpha) = \alpha +
\alpha^2 + \alpha^4$ for all $\alpha \in F$. It is not difficult to
see that $Tr$ is an additive homomorphism. Also, one can show that
$\ker {Tr} = \{ \alpha + \alpha^2 \mid \alpha \in F \}$.  Note that
$\card {\ker {Tr}} = 4$, and in fact, $F = Z_2 \oplus \ker {Tr}$. If
we fix an element $b \in F$, we can define a map $f : F \rightarrow
F$ by $f (a) = a b^4 + a^4 b$.  It is easy to show that $f$ is an
additive homomorphism whose kernel is $\{ 0, b \}$.

Consider the skew polynomial ring $F \{ X \}$ where we define the
multiplication of the indeterminant $X$ with elements of $F$ by $X
\alpha = \alpha^2 X$ for all $\alpha \in F$. We then define the ring
$R$ by $R = F\{X\}/(X^5)$, and we let $x$ be the image of $X$ in
$R$.  It is not difficult to see that the Jacobson radical $J$ of
$R$ has the form $\{ \alpha_1 x + \alpha_2 x^2 + \alpha_3 x^3 +
\alpha_4 x^4 \mid \alpha_i \in F \}$, and we let $Q = 1 + J$ as a
subgroup of the group of units of $R$. Observe that $|Q| = 2^{12}$.
We set $Q^i = 1 + J^i$ for $i = 1, \dots, 5$. We know from
\cite{coprime} that $Q^i$ is a subgroup of $Q$.  Also, $Q = Q^1$ and
$Q^5 = 1$.

Let $C = F^x$, the multiplicative group of $F$, and we see that $C$
embeds in the group of units of $R$ by $c \in C$ maps to the
constant polynomial $c 1$.  It follows that $C$ acts by conjugation
on $Q$ by $(1 + \alpha_1 x + \alpha_2 x^2 + \alpha_3 x^3 + \alpha_4
x^4)^c = 1 + \alpha_1 c x + \alpha_2 c^3 x^2 + \alpha_3 x^3 +
\alpha_4 c x^4$.  Observe that $C$ acts irreducibly on $Q/Q^2$,
$Q^2/Q^3$, and $Q^4/1$, but $C$ centralizes $Q^3/Q^4$.  Let ${\cal
G}$ be the Galois group of $F$ over $Z_2$, and we have an action of
${\cal G}$ on $QC$ by acting coordinate-wise on the elements of $Q$
and on $C$ as on $F$.  Notice that $C{\cal G}$ is a Frobenius group
of order $21$.

The group $Q$ has a $C{\cal G}$-invariant subgroup $P$ of index $2$.
To see this, observe that $Q^3/Q^4 = \cent {Q^3/Q^4}{{\cal G}}
\times [Q^3,{\cal G}]Q^4/Q^4$ via Fitting's lemma. We have $\cent
{Q^3/Q^4} {{\cal G}} = \{ 1 + \alpha x^3 + Q^4 \mid \alpha \in
Z_2\}$, and we set $X/Q^4 = [Q^3,{\cal G}]Q^4/Q^4 = \{ 1 + \alpha
x^3 + Q^4 \mid \alpha \in \ker {Tr} \}$. We note that $|Q^3:X| = 2$.
One can also show that $[Q,Q^2]Q^4/Q^4 = X/Q^4$ (this takes some
work), and so, $X = [Q,Q^2]Q^4$.  Also, we see that
$$
X = \{ 1 + \alpha_3 x^3 + \alpha_4 x^4 \mid \alpha_3 \in \ker {Tr},
\alpha_4 \in F \}.
$$

Now, $C$ stabilizes $X$, so $C$ acts on $Q^2/X$.  Note that $Q^2/X$
is central in $Q/X$, and thus, $Q^2/X$ is abelian.  By Fitting's
lemma, we have $Q^2/X = \cent {Q^2/X}C \times [Q^2,C]X/X$. We see
that $\cent {Q^2/X}C = Q^3/X$. Let $Y = [Q^2,C]X$, and observe that
$|Q^2:Y| = |Q^3:X| = 2$.  We see that $Q^2/Y = \cent {Q/Y}C$. Also,
we claim that
$$
Y = \{ 1 + \alpha_2 x^2 + \alpha_3 x^3 + \alpha_4 x^4 \mid \alpha_2,
\alpha_4 \in F, \alpha_3 \in \ker {Tr} \}.
$$

Using an argument similar to the one found in Theorem 4.5 of
\cite{family}, we can prove that $Y = Q'Q^4$.  In particular, $Q/Y$
is abelian, so $Q/Y = \cent {Q/Y}C \times [Q,C]Y/Y$.  We take $P =
[Q,C]Y$. Observe that $P' = Q' = Y$. We set $P^i = P \cap Q^i$, so
$P^1 = P = [Q,C]Y$, $P^2 = Y$, $P^3 = X$, and $P^4 = Q^4$. We have
$Z_1 (P) = P^4$, $Z_2 (P) = P^3$, and $Z_3 (P) = P^2$. In
particular, $P/Z_2 (P) = P/P^3 \cong Q/Q^3$ which is isomorphic to
the Suzuki $2$-group of type A and order $64$.  This proves that $P$
is an iterated Suzuki $2$-group.

Our group $G$ is the semidirect product of $C{\cal G}$ acting on
$P$.

\section{The theorem}

Following \cite{text}, we will say that $g \in G$ is a {\it real}
element if $g$ is conjugate to $g^{-1}$.  A conjugacy class of $G$
is called {\it real} if it consists of real elements.  In Problem
6.13 of \cite{text}, one is asked to prove that the number of real
irreducible characters of $G$ equals the number of real classes in
$G$.

\bigskip

\begin{theorem}
The group $G$ has exactly three real irreducible characters.
\end{theorem}

\begin{proof}
We now work to prove that $G$ has exactly $3$ irreducible
characters.  Now $G/P$ is the Frobenius group of order $21$, and as
such, the principal character is the only real valued character in
$\irr {G/P}$.  Let $N/P$ be the Frobenius kernel of $G/P$, and
observe that $|G:N| = 3$.  We note that $N/P^3$ is one of the
examples in \cite{threereal} where $\irr {N/P^3}$ has two real
characters.  It is not difficult to show that this will imply that
$\irr {G/P^3}$ has two real characters.  Working as in
\cite{family}, one can show that the nonprincipal characters in
$\irr {P^3/P^4}$ are fully ramified with respect to $P/P^3$, so the
characters in $\irr {P/P^4 \mid P^3/P^4}$ are real-valued.  If
$\theta \in \irr {P/P^4}$, then the stabilizer of $\theta$ in $G$ is
$N$.  Now, $\theta$ has a unique real extension to $N$.  Hence,
$\theta^G$ has a unique real irreducible constituent.  Note that
${\cal G}$ acts transitively on $\irr {P^3/P^4} - \{ 1 \}$, so,
$\irr {G/P^4 \mid P^3/P^4} = \irr {G \mid \theta}$. We conclude that
$\irr {G/P^4}$ has three real characters.

We will show that $G$ has $3$ real irreducible characters by showing
that $G$ has $3$ real classes.  Notice that $G/P^4$ has $3$ real
classes.  Also, it is not difficult to see that the nonidentity
elements of $P^2/P^4$ consist of involutions.  Notice that the
elements of $P^3/P^4 \setminus P^4/P^4$ and $P^2/P^4 \setminus
P^3/P^4$ must be in different conjugacy classes.  It follows that
all of the real elements of $G/P^4$ lie in $P^2/P^4$.  Since a real
element of $G$ gives rise to a real element of $G/P^4$, it follows
that the all real elements of $G$ lie in $P^2$.  In addition, notice
since $\gpindex GP$ is odd, that any element of $P^2$ that is real
in $G$ must in fact be real in $P$.

A typical element of $a \in Q$ has the form $a = 1 + a_1 x + a_2 x^2
+ a_3 x^3 + a_4 x^4$ with $a_i \in F$.  It is not difficult to see
that
$$
a^2 = 1 + a_1^3 x^2 + (a_1 a_2^2 + a_1^4 a_2) x^3 + (a_1 a_3^2 +
a_2^5 + a_1 a_3) x^4.
$$
It follows that $a^2 = 1$ if and only if $a_1 = a_2 = 0$.  Hence,
all of the involutions in $Q$ lie in $Q^3$, and all of the
nonidentity elements of $Q^3$ are involutions.  This implies that
all of the involutions in $P$ lie in $P^3$, and all of the
nonidentity elements of $P^3$ are involutions.  We know that $C$
acts transitively on $P^4 \setminus \{ 1 \}$.  So we have that $\{ 1
\}$ and $P^4 \setminus \{ 1 \}$ are two real conjugacy classes in
$G$.

Consider an element $b \in P^2$.  We know that $b = 1 + b_2 x^2 +
b_3 x^3 + b_4 x^4$ with $b_2, b_4 \in F$ and $b_3 \in \ker {Tr}$.
Let $a$ be any element of $Q$ as in the previous paragraph.
Computing, we see that
$$
ab = 1 + a_1 x + (a_2 + b_2)x^2 + (a_3 + b_3 + a_1 b_2^2)x^3 + (a_4
+ b_4 + a_1b_3^2 + a_2b_2^4)x^4.
$$
and
$$
ba = 1 + a_1 x + (a_2 + b_2)x^2 + (a_3 + b_3 + b_2 a_1^4)x^3 + (a_4
+ b_4 + b_3a_1 + b_2a_2^4)x^4.
$$
It follows that $ab = ba$ if and only if $(1)~a_1 b_2^2 = b_2 a_1^4$
and $(2)~a_1 b_3^2 + a_2 b_2^4 = b_3 a_1 + b_2 a_2^4$.  If $b \in
P^3 \setminus P^4$, then $b_2 = 0$ and $b_3 \ne 0$.  We see that $ab
= ba$ if and only if $a_1 b_3^2 = b_3 a_1$.  Recall that $b_3 \in
\ker {Tr}$, so $b_3 \ne 1$.  It follows that $b_3 \ne b_3^3$, and
hence, $a_1 = 0$.  We conclude $\cent Qb = Q^2$, and so, $\cent Pb =
P^2$. We then have that $\cent Gb = P^2C$.  It follows that
$\gpindex G{\cent Gb} = 8 \cdot 3 = 24$, and so the conjugacy class
containing $b$ has size $24$.  On the other hand, $\card {P^3
\setminus P^4} = \card {P^3} - \card {P^4} = 32 - 8 = 24$.  Thus,
$P^3 \setminus P^2$ form the third real conjugacy class of $G$.

We will now show that these are the only real conjugacy classes in
$G$.  As we have seen, $P^2$ contains all of the real classes in
$G$.  Thus, we consider $b \in P^2 \setminus P^3$.  This implies
that $b_2 \ne 0$.  It is not difficult to see that the only possible
solutions for $a_1$ in the equation $a_1 b_2^2 = a_1^4 b_2$ are $a_1
= 0$ and $a_1 = b_2^5$.  If $a_1 = 0$, then equation $(2)$ becomes
$a_2 b_2^4 = a_2^4 b_2$.  It is not difficult to see that the only
solutions for $a_2$ in this equation are $a_2 = 0$ and $a_2 = b_2$.
On the other hand, if $a_1 = b_2^5$, then equation $(2)$ becomes
$b_2^5 b_3^2 + a_2 b_2^4 = b_3 b_2^5 + b_2 a_2^4$.  This implies
that $a_2 b_2^4 + a_2^4 b_2 = b_2^5 (b_3 + b_3^2)$.  We have seen
that the map $a_2 \mapsto a_2 b_2^4 + a_2^4 b_2$ is an additive
homomorphism with a kernel of size $2$.  This implies that our
equation has $2$ solutions for $a_2$.  It is not difficult to see
that $P^3 \le \cent Pb$.  It follows that $\gpindex P{\cent Pb} =
16$.  Since $\cent Pb$ will contain a conjugate of ${\cal G}$, we
deduce that $\gpindex G{\cent Gb} = 16 \cdot 7 = 112$.

Note that $\card {P^2 \setminus P^3} = \card {P^2} - \card {P^3} =
256 - 32 = 224$.  Therefore, $P^2 \setminus P^3$ consists of $2$
conjugacy classes of $G$.  This implies that either every element of
$P^2 \setminus P^3$ is conjugate to its inverse, or none of the
elements in $P^2 \setminus P^3$ are conjugate to their inverse.
Thus, we will fix $b = 1 + x^2$.  It suffices to show that $b$ is
not conjugate to $b^{-1}$ in $G$.   As we have seen, it suffices to
show that $b$ is not conjugate to $b^{-1}$ in $P$.  We will show
that in fact $b$ is not conjugate to $b^{-1}$ in $P$.

By way of contradiction, we suppose that $b$ is conjugate to
$b^{-1}$ in $Q$.  Hence, there is an element $a \in Q$ so $b^a =
b^{-1}$.  Observe that $b^{-1} = 1 + x^2 + x^4$. As in the previous
paragraph, we can write $a = 1 + a_1 x + a_2 x^2 + a_3 x^3 + a_4
x^4$. Let $c = b^a$, and observe that $ac = ba$. Observe that $ba =
1 + a_1 x + (a_2 + 1) x^2 + (a_3 + a_1^4) x^3 + (a_4 + a_2^4) x^4$.
We know that $c \in P^2$, so $c = 1 + c_2 x^2 + c_3 x^3 + c_4 x^4$,
and
$$
ac = 1 + a_1 x + (a_2 + c_2) x^2 + (a_3 + c_3 + a_1 c_2^2) x^3 +
(a_4 + c_4 + a_1 c_3^2 + a_2 c_2^4) x^4.
$$
This implies that $c_2 = 1$ and $c_3 = a_1^4 + a_1$ and
$$
c_4 = a_2^4 + a_1 (a_1^4 + a_1)^2 + a_2 = (a_2^4 + a_2) + (a_1^4 +
a_1)^2.
$$
If $c = b^{-1}$, then matching coefficients, we must have $a_1^4 +
a_1 = 0$ and $(a_2^4 + a_2) + (a_1^4 + a_1)^2 = 1$.  Combining these
two equations, we obtain $a_2^4 + a_2 = 1$.  However, $a_2^4 + a_2
\in \ker {Tr}$ and $1 \not\in \ker {Tr}$, and so, we have a
contradiction.  This contradiction implies that $b$ is not conjugate
to its inverse.  Therefore, $G$ has $3$ real classes and hence, $3$
real valued irreducible characters.
\end{proof}

A natural extension of this group is to start our construction by
starting with the ring $F\{X\}/\langle X^6 \rangle$, and then
completing the construction as before.  To see that the resulting
group has more than three real irreducible characters, we note that
$$
\{ 1 + a_3 x^3 + a_4 x^4 + a_5 x^5 \mid a_3 \in \ker {Tr}, a_4, a_5
\in F \}
$$
is a characteristic subgroup of the resulting group.  It
is not difficult to see that this subgroup has all its nonidentity
elements are involutions.  Also, it is not difficult to see that the
involutions form at least $3$ conjugacy classes, and so, the group
has at least $4$ real conjugacy classes.

\section{Appendix: Magma code}

\noindent Next, include code (in Magma) to generate $G$ as PC-group:

\medskip

$>g :=$ PolycyclicGroup $<s,c,x1,x2,x3,y1,y2,y3,w1,w2,z1,z2,z3 |$

$>   s^\wedge3,c^\wedge7,c^\wedge s=c^\wedge4,$

$>   x1^\wedge2=y1,x2^\wedge2=y2,x3^\wedge2=y3,y1^\wedge2=z1,
y2^\wedge2=z2,y3^\wedge2=z3,$

$>   w1^\wedge2,w2^\wedge2,z1^\wedge2,z2^\wedge2,z3^\wedge2,$

$>   x2^\wedge x1=x2*y2*y3*w1*w2*z1*z2*z3,x3^\wedge
x1=x3*y1*y2*w1*z2,$

$>   x3^\wedge x2=x3*y1*w1*w2*z1*z3,$

$>   y2^\wedge x1=y2*w1,y3^\wedge x1=y3*w2,$

$>   y2^\wedge y1=y2*z3,y3^\wedge y1=y3*z2*z3,y3^\wedge
y2=y3*z1*z2,$

$>   y1^\wedge x2=y1*w1*w2*z1*z2*z3,y3^\wedge x2=y3*w1*z1,$

$>   y1^\wedge x3=y1*w1*z1*z2,y2^\wedge x3=y2*w1*w2*z2,$

$>   w1^\wedge x1=w1*z3,w2^\wedge x1=w2*z2*z3,w1^\wedge
x2=w1*z1*z2,w2^\wedge x2=w2*z1*z2*z3,$

$>   w1^\wedge x3=w1*z2*z3,w2^\wedge x3=w2*z1*z3,$

$>   x1^\wedge c=x1*x2*x3*w1*z1,x2^\wedge
c=x1*x3*y1*y3*w1*z3,x3^\wedge c=x1,$

$>   y1^\wedge c=y1*y2*z1*z2,y2^\wedge c=y2*y3*z2*z3,y3^\wedge
c=y1,$

$>   w1^\wedge c=w1*z2,w2^\wedge c=w2*z1*z2*z3,$

$>   z1^\wedge c=z1*z2*z3,z2^\wedge c=z1*z3,z3^\wedge c=z1,$

$>   x1^\wedge s=x1*z3,x2^\wedge s=x2*x3*y2*y3*w1*w2*z1*z3,x3^\wedge
s=x2,$

$> y2^\wedge s=y1*y2*y3,y3^\wedge s=y2,w1^\wedge s=w1*w2,w2^\wedge
s=w1,z2^\wedge s=z2*z3,z3^\wedge s=z2 >;$

\medskip

Using Magma, we computed the character table of the above group.
This group has a character table with $33$ irreducible characters.
Given this size, we have not tried to reproduce this character table
in this note. The character degrees of $G$ are $\{ 1, 3, 7, 14, 24,
56 \}$. There are three real irreducible characters. One each of
degrees $1$, $7$, and $24$.


\end{document}